\documentclass[12pt]{amsart}
\usepackage{verbatim}
\usepackage{amsmath}
\usepackage{amsthm}
\usepackage{amssymb}
\usepackage{enumerate}

\newif\ifdeveloping

\newtheorem{theorem}{Theorem}
\newtheorem{proposition}[theorem]{Proposition}

\newtheorem{corollary}[theorem]{Corollary}

\newtheorem{problem}[theorem]{Problem}

\newtheorem{qtheorem}{Theorem}

\theoremstyle{definition}

\newtheorem{definition}[theorem]{Definition}

\theoremstyle{remark}
  
\newtheorem{example}[theorem]{Example}{}

\ifdeveloping
\usepackage[notref,notcite]{showkeys}
\fi

\newcommand{\prtime}{{\count0=\time\divide\count0 by 60
\count1=-\count0\multiply\count1 by 60 \advance\count1 by \time
\the\count0:\the\count1} }

\def\myheads#1;#2;{
\pagestyle{myheadings} \markboth{{\sc\hfill
#1\hfill\protect\makebox[0cm][r]{\rm\today; \prtime}}}
{{\sc\protect\makebox[0cm][l]{\rm\today;\ \prtime}\hfill #2\hfill}}
\thispagestyle{myheadings} }

\newcommand{\subs}{\subset}

\def\<{\left\langle}
\def\>{\right\rangle}

\def\br#1;#2;{\bigl[ {#1} \bigr]^ {#2} }
\def\bc#1;#2;{\bigl( {#1} \bigr)^ {#2} }

\newcommand{\pir}{$\pi$-irreducible}

\newcommand{\p}{p\,\pi\chi}

\def\to{\longrightarrow}

\def\Int{\operatorname{Int}}

\theoremstyle{plain}

\begin{document}
\author[I. Juh\'asz]{Istv\'an Juh\'asz}
\address{Alfr{\'e}d R{\'e}nyi Institute of Mathematics}
\email{juhasz@renyi.hu}

\author[Z. Szentmikl\'ossy]{Zolt\'an Szentmikl\'ossy}
\address{E\"otv\"os Lor\'ant University, Department of  Analysis,
1117 Budapest,
 P\'azm\'any P\'eter s\'et\'any 1/A,
 Hungary}
\email{zoli@renyi.hu}
\thanks{Research supported
by OTKA grant no. 61600.}



\title[Projective $\pi$-character]{Projective $\pi$-character bounds the order of a $\pi$-base}

\begin{abstract}
All spaces below are Tychonov. We define the projective
$\pi$-character $\p(X)$ of a space $X$ as the supremum of the values
$\pi\chi(Y)$ where $Y$ ranges over all continuous images of $X$. Our
main result says that every space $X$ has a $\pi$-base whose order
is $\le \p(X)$, that is every point in $X$ is contained in at most
$\p(X)$-many members of the $\pi$-base. Since $\p(X) \le t(X)$ for
compact $X$, this provides a significant generalization of a
celebrated result of Shapirovskii.

\end{abstract}

\maketitle

Arhangel'skii has recently introduced in \cite{A} the concept of a
space of {\em countable projective $\pi$-character} and noticed that
any compact space of countable tightness has countable projective
$\pi$-character. Then he showed that a compact space of countable
projective $\pi$-character having $\omega_1$ as a caliber is
separable, thereby strengthening Shapirovskii's analogous result for
countably tight compacta. Note that Shapirovskii's theorem is a
trivial corollary of his more general result establishing that any
countably tight compactum has a point-countable $\pi$-base, or more
generally: any compactum $X$ has a $\pi$-base of order at most
$t(X)$, see \cite{J} or \cite{Sh}.

In this paper we use the general concept of projective
$\pi$-character to give the following significant generalization of
this stronger result of Shapirovskii: Any Tychonov space has a
$\pi$-base of order at most the projective $\pi$-character of the
space. Not only is this result stronger for compacta, by replacing
tightness with projective $\pi$-character that is smaller, but
somewhat surprisingly it extends to all Tychonov spaces.

Let $\varphi$ be any cardinal function defined on a class
$\mathcal{C}$ of topological spaces. We define the projective
version $p\,\varphi$ of $\varphi$ as follows. For any $X \in
\mathcal{C}$ we let $p\,\varphi(X)$ be the the supremum of the
values $\varphi(Y)$ where $Y$ ranges over all continuous images of
$X$ in $\mathcal{C}$. In particular, we shall consider the case in
which $\varphi = \pi \chi$, the $\pi$-character defined on the class
of Tychonov spaces. It is easy to show that then a Tychonov space
$X$ has countable projective $\pi$-character in the sense of
\cite{A} iff $\p(X) \le \omega$.

Also, as was already mentioned before, if $X$ is compact Hausdorff
then we have $\p(X) \le t(X).$ In fact, this follows because $t(Y)
\le t(X)$ for any continuous image of $X$ and $\pi\chi(Y) \le t(Y)$
for every compact Hausdorff  $Y$. But are $\pi\chi(X)$ and $t(X)$
the same? Arhangel'skii asked, more specifically, if there is a
compactum of countable projective $\pi$-character that is not
countably tight, see \cite{A}, problem 7. The next example yields
such a compactum.

\begin{example}
Let $X$ be a compactification of $\omega$ whose remainder is
(homeomorphic to) $\omega_1+1.$ Then $\p(X) \le \omega < t(X).$
\end{example}

\begin{proof}
It is obvious that $t(\omega_1,X) = t(X) = \omega_1$. To see $\p(X)
\le \omega$, consider any continuous surjection $f : X \to Y$. If
$f(\omega_1) = p$ is an isolated point in $Y$ then there is an
$\alpha < \omega_1$ such that $f$ is constant on the interval
$[\alpha,\omega_1]$, hence $Y$ is countable and compact, so,
trivially, $\pi\chi(Y) \le w(Y) = \omega$.

If, however, $p$ is not isolated then $Y$ has a countable dense
subset $S$ with $p \notin S$. Then there is a closed $G_\delta$ set
$F$ such that $p \in F \subs Y \backslash S$ and again we can find
an $\alpha < \omega_1$ such that $f[\alpha,\omega_1] = F$. But then
$G = Y \backslash F$ is countable and dense open in $Y$, moreover
$w(G) = \omega$ because every countable locally compact space is
second countable. So we have $\pi\chi(Y) \le \pi(Y) = w(G) =
\omega$.

\end{proof}

We recall from \cite{J} that $\pi sw(X)$ denotes the {\em
$\pi$-separating weight} of a space $X$, that is the minimum order
of a $\pi$-base of $X$, see p. 74 of \cite{J}.

With this we may now formulate our main result as follows.

\begin{theorem}\label{main}
For any Tychonov space $X$ we have $\pi sw(X) \le \p(X).$ In
particular, any Tychonov space of countable projective
$\pi$-character has a point-countable $\pi$-base.
\end{theorem}

Our proof of theorem \ref{main} will go along similar lines as
Shapirovskii's proof of the weaker result $\pi sw(X) \le t(X)$ for
compact spaces, however the role of irreducible maps in it will be
played by a new type of maps that we shall call $\pi$-irreducible.
So we shall first define and deal with these maps.

\begin{definition}
Let $f$ be a continuous map of $X$ {\em onto $Y$}. We say that the
map $f$ is \pir\ if for every proper closed subset $F \subs X$ its
image $f[F]$ is not dense in $Y$.
\end{definition}

Clearly, an onto map $f$ is \pir\ iff the $f$-image of a non-dense
set is non-dense. Also, it is obvious that a closed map is \pir\
iff it is irreducible, consequently the two concepts coincide for
maps between compact Hausdorff spaces.

The following proposition will be used in the proof of theorem
\ref{main} and explains our terminology.

\begin{proposition}\label{pir}
Let $f$ be a continuous map of $X$ onto $Y$. Then the following five
statements (1)--(5) are equivalent.
\begin{enumerate}
\item $f$ is \pir\ ;
\smallskip
\item for every $\pi$-base $\mathcal{B}$ of $X$ and for every $B \in
\mathcal{B}$ the $f$-image of its complement, $f[X \backslash B]$,
is not dense in $Y$;
\smallskip
\item there is a $\pi$-base $\mathcal{B}$ of $X$ such that for every $B \in
\mathcal{B}$ the $f$-image $f[X \backslash B]$ is not dense in $Y$;
\smallskip
\item for every $\pi$-base $\mathcal{C}$ of $Y$ the family $\{f^{-1}(C) : C \in
\mathcal{C}\}$ is a $\pi$-base of $X$ ;
\smallskip
\item there is a $\pi$-base $\mathcal{C}$ of $Y$ such that $\{f^{-1}(C) : C \in
\mathcal{C}\}$ is a $\pi$-base of $X$ .

\end{enumerate}

\end{proposition}

\begin{proof}
We shall show (3)$\Rightarrow$(4) and (5)$\Rightarrow$(1) only
because the other three implications of the cycle are trivial.

So, let $\mathcal{B}$ be as in (3) and $\mathcal{C}$ be any
$\pi$-base of $Y$. For every non-empty open set $U$ in $X$ choose $B
\in \mathcal{B}$ with $B \subs U$. Then there is a $C \in
\mathcal{C}$ such that $C \cap f[X \backslash B] = \emptyset$, and
hence $f^{-1}(C) \subs B \subs U$.

Now, let $\mathcal{C}$ be as in (5) and $F$ be a proper closed
subset of $X$. Then there is a $C \in \mathcal{C}$ with $F \cap
f^{-1}(C) = \emptyset$, consequently we have $f[F] \cap C =
\emptyset\,$ and so $f[F]$ is not dense in $Y$.
\end{proof}

\begin{corollary}\label{pie}
If $f : X \to Y$ is \pir\ then $\,\,\pi(X) = \pi(Y).$
\end{corollary}

\begin{proof}
$\pi(X) \le \pi(Y)$ is immediate from part (4) of proposition
\ref{pir}. To see $\pi(X) \ge \pi(Y)$ first note that for any
non-empty open $U \subs X$ the interior of $f[U]$ in $Y$ is
non-empty. So for any $\pi$-base $\mathcal{B}$ of $X$ the family $\{
\Int_Y(f[B]) : B \in \mathcal{B} \}$ is a $\pi$-base of $Y$. Indeed,
this is because if $V$ is non-empty open in $Y$ and $B \in
\mathcal{B}$ with $B \subs f^{-1}(V)$ then $f[B] \subs V.$
\end{proof}

We now consider another key ingredient of the proof of our main
result: certain specially embedded subspaces of Tychonov cubes. As
usual, we shall denote the unit interval $[0,1]$ by $I$. The members
of the Tychonov cube $I^\kappa$ will be construed as functions from
$\kappa$ to $I$. So if $x \in I^\kappa$ and $\alpha < \kappa$ then
$x\upharpoonright\alpha$ is the projection of $x$ to the subproduct
$I^\alpha$.

\begin{definition}
We say that $Y \subs I^\kappa$ is 0-embedded in the Tychonov cube
$I^\kappa$ if $$\{y\upharpoonright\alpha : y \in Y \mbox{ and }
y(\alpha) = 0 \}$$ is dense in the projection
$Y\upharpoonright\alpha = \{ y\upharpoonright\alpha : y \in Y \}$
for every $\alpha < \kappa$.
\end{definition}

We now present two results concerning 0-embedded subspaces of
Tychonov cubes which will be crucial in the proof of our main
theorem and are also interesting in themselves.

\begin{theorem}\label{0-em}
Assume that $Y$ is 0-embedded in the Tychonov cube $I^\kappa$ where
$\kappa$ is a regular cardinal and $y \in Y$ is such that $y(\alpha)
> 0$ for all $\alpha < \kappa$. Then $\pi\chi(y,Y) = \kappa.$
\end{theorem}

\begin{proof}
Of course, only $\pi\chi(y,Y) \ge \kappa$ needs to be proven. To see
this, let $\mathcal{U}$ be any family of elementary open sets in
$I^\kappa$ such that $|\mathcal{U}| < \kappa$ and $U \cap Y \ne
\emptyset$ for all $U \in \mathcal{U}$. Every elementary open set $U
\in \mathcal{U}$ is supported by a finite subset of $\kappa$, hence
the regularity of $\kappa$ implies the existence of an ordinal
$\alpha < \kappa$ such that the support of each $U \in \mathcal{U}$
is included in $\alpha$.

Since $Y$ is 0-embedded in $I^\kappa$, this implies that for every
$U \in \mathcal{U}$ we may pick a point $y_U \in U \cap Y$ such that
$y_U(\alpha) = 0$. But then $y(\alpha) > 0$ clearly implies that the
point $y$ is not in the closure of the set $\{ y_U : U \in
\mathcal{U} \}$, consequently $\mathcal{U}$ cannot be a local
$\pi$-base at $y$ in $Y$, completing the proof.
\end{proof}

From theorem \ref{0-em} we can immediately obtain the following
useful corollary about the projective $\pi$-character of
0-embedded subspaces of Tychonov cubes.

\begin{corollary}\label{0-emc}
If $Y$ is 0-embedded in the Tychonov cube $I^\kappa$ then for every
$y \in Y$ we have
$$\p(Y) \ge \big|\{ \alpha : y(\alpha) > 0 \}\big|.$$
\end{corollary}

Our next result shows that every Tychonov space admits a
$\pi$-irrredu-cible map onto a suitable 0-embedded subspace of a
Tychonov cube.

\begin{theorem}\label{pi0}
Let $X$ be any Tychonov space of $\pi$-weight $\pi(X) = \kappa$.
Then there is a \pir\ map $f$ of $X$ onto a 0-embedded subspace $Y$
of the Tychonov cube $I^\kappa$.
\end{theorem}

\begin{proof}
To begin with, let us choose a $\pi$-base $\mathcal{B}$ of $X$ with
$|\mathcal{B}| = \kappa$ and fix a well-ordering $\prec$ of
$\mathcal{B}$ of order-type $\kappa$.

We shall define by transfinite induction on $\alpha < \kappa$ the
co-ordinate maps $g_\alpha = p_\alpha \circ f : X \to I$, where
$p_\alpha(y) = y(\alpha)$ is the $\alpha$th co-ordinate projection,
and sets $B_\alpha \in \mathcal{B}$. So assume that $\alpha <
\kappa$ and for all $\beta < \alpha$ the maps $g_\beta : X \to I$
and the sets $B_\beta \in \mathcal{B}$ have been defined.

Let $f_\alpha : X \to I^\alpha$ be the map whose $\beta$th
co-ordinate map is $g_\beta$ for all $\beta < \alpha$ and set
$Y_\alpha = f_\alpha[X]$. Then, in view of corollary \ref{pie},
the map $f_\alpha : X \to Y_\alpha$ cannot be \pir\ because
$\pi(Y_\alpha) < \kappa = \pi(X)$, hence using part (2) of
proposition \ref{pir} there is a member $B \in \mathcal{B}$ for
which $f_\alpha[X \backslash B]$ is dense in $Y_\alpha$. Let
$B_\alpha$ be the $\prec$-first such member of $\mathcal{B}$. We
then define $g_\alpha : X \to I$ as any continuous function that
is identically $0$ on $X \backslash B_\alpha$ and takes the value
$1$ at some point in $B_\alpha$. As was intended, with the
induction completed we let $f : X \to I^\kappa$ be the unique map
having the $g_\alpha$ for $\alpha < \kappa$ as its co-ordinate
functions and we also set $Y = f[X]$.

Note first that if $\beta < \alpha$ then $B_\beta \prec B_\alpha$.
Indeed, since we have $Y_\beta = Y_\alpha\upharpoonright\beta$, the
density of $f_\alpha[X \backslash B_\alpha]$ in $Y_\alpha$ implies
that $f_\beta[X \backslash B_\alpha]$ is dense in $Y_\beta$,  hence
$B_\alpha \prec B_\beta$ would contradict the choice of $B_\beta$.
Moreover, by our construction, $f_{\beta+1}[X \backslash B_\beta]$
is not dense in $Y_{\beta+1}$ and consequently $f_\alpha[X
\backslash B_\beta]$ is not dense in $Y_\alpha$, which implies
$B_\alpha \ne B_\beta$.

Since $\mathcal{B}$ is of order type $\kappa$ under $\prec$, it
follows from this that for every $B \in \mathcal{B}$ there is an
$\alpha < \kappa$ with $B \prec B_\alpha$. But then, by the choice
of $B_\alpha$ we have that $f_\alpha[X \backslash B]$ is not dense
in $Y_\alpha = Y\upharpoonright\alpha$ and hence $f[X \backslash
B]$ cannot be dense in $Y$. Using part (3) of proposition
\ref{pir} this implies that $f$ is indeed a \pir\ map of $X$ onto
$Y$.

Finally, by our construction, for every $\alpha < \kappa$ the image
$f_\alpha[X \backslash B_\alpha]$ is dense in $Y_\alpha = Y
\upharpoonright \alpha$, moreover we have
$$f_\alpha[X \backslash B_\alpha] \subs \{y\upharpoonright\alpha : y
\in Y \mbox{ and } y(\alpha) = 0 \},$$ consequently $Y$ is indeed
0-embedded in $I^\kappa$.
\end{proof}

Let us now recall that the $\kappa$th $\Sigma_\lambda$-power of $I$,
denoted by $\Sigma_\lambda(I,\kappa)$, is the subspace of $I^\kappa$
consisting of all points whose support is of size at most $\lambda$.
The support of a point $y \in I^\kappa$ is the set $\{ \alpha <
\kappa : y(\alpha) > 0 \}$. Thus, from theorem \ref{pi0} and from
corollary \ref{0-emc}, moreover from the trivial fact that $\p(Y)
\le \p(X)$ if $Y$ is any continuous image of $X$, we immediately
obtain the following result.

\begin{corollary}\label{sig}
If $X$ is a Tychonov space such that $\pi(X) = \kappa$ and $\p(X) =
\lambda$ then some \pir\ image $Y$ of $X$ embeds into
$\Sigma_\lambda(I,\kappa)$.
\end{corollary}

This corollary is clearly a strengthening of the following result of
Shapirovskii from \cite{Sh} (see also 3.22 of \cite{J}) : If $X$ is
compact Hausdorff then some irreducible image of $X$ embeds into a
$\Sigma_{t(X)}$-power of $I$.

The proof of our main theorem \ref{main} can now be easily
established by recalling the following result of Shapirovskii from
\cite{Sh} (see also\cite{J}, 3.24).

\begin{qtheorem}[Shapirovskii]
If the space $Y$ embeds into a $\Sigma_\lambda$-power of $I$ then
$\pi sw(Y) \le \lambda$.
\end{qtheorem}

\begin{proof}[\bf Proof of theorem \ref{main}]
Now, to prove theorem \ref{main}, consider any Tychonov space $X$.
By corollary \ref{sig} then there is a \pir\ map $f$ of  $X$ onto a
space $Y$ such that $Y$ embeds into a $\Sigma_\lambda$-power of $I$,
where $\lambda = \p(X)$. By the previous theorem of Shapirovskii,
the space $Y$ has a $\pi$-base $\mathcal{C}$ of order at most
$\lambda$. But by part (4) of proposition \ref{pir}, then the family
$\{f^{-1}(C) : C \in \mathcal{C}\}$ is a $\pi$-base of $X$ that
clearly has the same order as $\mathcal{C}$.
\end{proof}

The following result is then an immediate consequence of theorem
\ref{main}.

\begin{corollary}\label{cal}
Let $X$ be any Tychonov space and $\kappa > \p(X)$ be a cardinal
such that $\kappa$ is a caliber of $X$. Then $\pi(X) < \kappa$.
\end{corollary}

Since $t(X) \ge \p(X)$ for a compact Hausdorff space $X$, this
corollary implies Shapirovskii's theorem which says that if $t(X)^+$
is a caliber of such a space $X$ then $\pi(X) \le t(X)$. Moreover,
it also extends from compacta to all Tychonov spaces Arhangel'skii's
result from \cite{A} saying that spaces of countable projective
$\pi$-character and having $\omega_1$ as a caliber are separable.

Let us conclude this paper by pointing out that neither theorem
\ref{main} nor corollary \ref{cal} remain valid if the projective
$\pi$-character $\p$ is replaced by simple $\pi$-character
$\pi\chi$ in them. In fact, it has recently been shown in
\cite{JSSz} that there are even first countable spaces whose
$\pi$-separating weight is as large as you wish. Moreover, in the
same paper it was also shown that it is consistent to have first
countable spaces with caliber $\omega_1$ which have uncountable
$\pi$-weight (or equivalently, density). However, since first
countability implies countable tightness, none of these examples
are (or could be) compact, so the following intriguing questions
remain open.

\begin{problem}
Let $X$ be a compact Hausdorff space of countable $\pi$-character.
Does $X$ have a point-countable $\pi$-base? If, in addition,
$\omega_1$ is a caliber of $X$, is then $X$ separable?
\end{problem}

\end{document}